\documentclass[11pt]{amsart}
\usepackage{graphicx,multirow}
\usepackage[mathcal]{euscript}
\usepackage{float}
\usepackage{amsmath,amsthm,amssymb,amsfonts,amscd,epsfig,latexsym,graphicx,textcomp}

\restylefloat{figure}

\theoremstyle{definition}

\theoremstyle{remark}

\numberwithin{equation}{section}

\newcommand{\BR}{\mathbb{R}}

\newcommand{\euler}{\mathbf{e}}

\begin{document}

\title[Small cube diagrams]{Small examples of Cube diagrams of knots}

\author[S. Baldridge]{Scott Baldridge}
\author[B. McCarty]{Ben McCarty}

\thanks{S. Baldridge was partially supported by  NSF Grant DMS-0748636.}

\address{Department of Mathematics, Louisiana State University \newline
\hspace*{.375in} Baton Rouge, LA 70817, USA} \email{\rm{sbaldrid@math.lsu.edu}}

\address{Department of Mathematics, Louisiana State University \newline
\hspace*{.375in} Baton Rouge, LA 70817, USA} \email{\rm{benm@math.lsu.edu}}

\date{July 30, 2009}

\begin{abstract}
In this short note we highlight some of the differences between cube diagrams and grid diagrams.  We also list examples of small cube diagrams for all knots up to 7 crossings and give some examples of links.

\end{abstract}

\maketitle

\bigskip
\section{Introduction}
\bigskip

Cube diagrams, introduced in Baldridge and Lowrance~\cite{Adam}, are 3-dimensional representations of knots or links.  We define cube diagrams carefully in Section 3.  The easiest way to imagine a cube diagram is to think of an embedding of a knot or link in a $[0,n]\times [0,n]\times [0,n]$ cube (using $xyz$ coordinates) for some positive integer $n$ such that the knot projection of the cube to each axis plane ($x=0$, $y=0$, and $z=0$) is a grid diagram.

\begin{figure}[H]

\includegraphics[scale = .4]{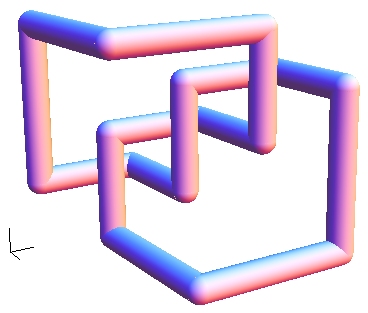}
\caption{A cube diagram for the Trefoil.}
\label{Intro:Trefoil}
\end{figure}

The integer $n$ is called the {\em size} of the cube diagram and we will refer to a cube diagram as ``small'' if its size is close to the crossing number of the knot.  Small cube diagrams are very useful for computing and testing for knot invariants using computers---calculating certain invariants of cube diagrams with size $n>20$ can be computationally intractable for most computers.  In this paper we discuss why small cube diagrams can be elusive to find for a given knot or link and we give examples of small cube diagrams ($n<20$) for knots with seven crossings or less.  In particular, the Appendix gives pictures of the cube diagrams we have generated.  A Mathematica program \cite{BLcode} can be downloaded from \verb+http://cubeknots.googlecode.com+ that can be used to calculate invariants of cube diagrams and also to rotate in 3-dimensions the examples in this paper.

\medskip

Grid diagrams and cube diagrams are useful representations in knot theory.  In 1996 Cromwell used grid diagrams as ways  to represent embeddings of knots in open books (c.f. \cite{Cromwell}).  He described a set of elementary grid moves that preserve topological knot type that can be used to check for knot invariants.  More recently grid diagrams were used in constructing a combinatorial version of knot Floer homology (c.f. \cite{knotfloer} and \cite{linkfloer}) and they provide a natural presentation of the front projection of a Legendrian knot (c.f. \cite{legend} and \cite{Ng}).  Cube diagrams share many similarities with grid diagrams.  Like grid diagrams, there is set of elementary cube moves that preserve the knot type and there is a combinatorial knot Floer homology that can be computed from a cube diagram \cite{Adam}.  However,  cube diagrams have also appeared independent of grid diagrams in the study of unstable Vassiliev theory \cite{Giusti}.

\bigskip
\section{Finding small cube diagrams: the lifting problem}
\bigskip

In order to find small cube diagrams, it is important to understand how they are different from grid diagrams.  In this section we discuss the problem of lifting a grid diagram to a cube diagram.

\medskip

A grid diagram $G$ is an $n \times n$ square grid decorated with $X$ and $O$ markings in such a way that every row (resp. column) contains exactly one $X$ and one $O$ marking.  To get an oriented knot or link projection from a grid diagram one draws edges from $X$ to $O$ in each column and from $O$ to $X$ in each row, taking the vertical segment as the over crossing at any intersection (cf.~\cite{Adam,Cromwell}).

\begin{figure}[H]
\includegraphics[scale = .1]{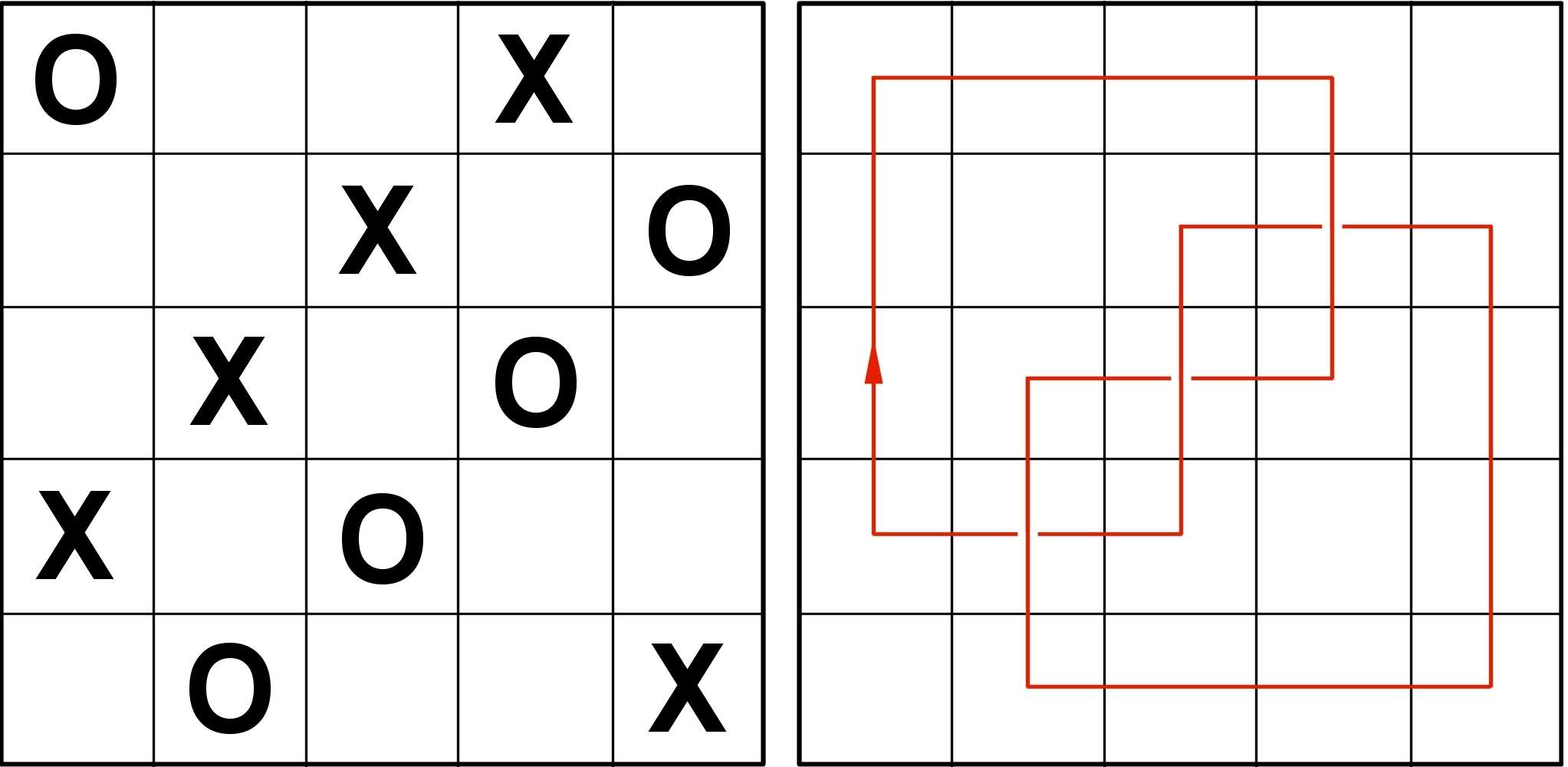}
\caption{Grid diagram:  $X=\{2,3,4,5,1\}$ and $O=\{5,1,2,3,4\}$.}
\label{Methodology:XO}
\end{figure}

Since each of the three projections of a cube diagram is a grid diagram, it is natural to think of a cube diagram as a lift of a grid diagram of one of the three projections.  Such a lift is clearly not unique.  What is not obvious is that many times such lifts do not exist: a given grid diagram does not necessarily lift to an embedding of a knot or link in a lattice in an $n\times n \times n$ cube such that the projection to each plane is a well-defined knot projection (it can lift to an embedding, but the intersections in the other two projections are not necessarily isolated double points).  For an example of a grid diagram that can not be the projection of such a lattice knot see \cite{Adam}.


\medskip

The lifting problem to a true cube diagram is harder than just finding a lift of a grid diagram in 3-space such that the other two projections are valid knot projections.  Even when a grid diagram does lift to such an embedding, that embedding is rarely a cube diagram and the chances that a grid diagram lifts to a cube diagram appears to decrease as the grid size increases.  To support this assertion we wrote a brute-force program to get statistics on grid diagrams that lift to cube diagrams.  The program looked at all size $n=5,6,7,8$ grid diagrams and tested whether each grid diagram was a nontrivial knot (we did not consider links) and whether each grid diagram lifted to a cube diagram.  The data for grid diagrams for nontrivial knots are presented below.

\medskip

\begin{center}
\begin{tabular}{|c|r|r|r|}
\hline
Grid & Total grids   & Number that & Percent\\
Size & of nontrivial & lift to cube & \mbox{}\\
\mbox{}& knots         & diagrams     & \mbox{}\\ \hline
5 & 10 & 3 & 30\%\\ \hline
6 & 972 & 261 & 27\% \\ \hline
7 & 85,022 & 19,722 & 23\% \\ \hline
8 & 8,077,072  & 1,589,447  & 19.7\% \\ \hline
\end{tabular}
\end{center}

\medskip


Note that creating grid diagrams with sizes above 8 is time intensive and the computation to rule out unknot grid diagrams grows quickly beyond the capabilities of most computers.  For example, the total number of grid diagrams of knots of size 8 including unknots is 101,606,400.  It is interesting to note that unknot grid diagrams do lift to cube diagrams more often:  of those 101 million grid diagrams, 72,109,568 of them or 71 percent lift to cube diagrams.  But that fact is because many unknot grid diagrams do not have any crossings and such grid diagrams often lift to cube diagrams.  Furthermore, the percentage of unknot grid diagrams that lift to cube diagrams also decreases as the grid size increases in the examples we have calculated.  To get the data for size 8 grid diagrams, we ran 5 computers simultaneously night and day for 1 week.   See the program \cite{BLcode} for other methods of measuring the probability of finding a cube diagram when the size is greater than 8.

\medskip

Therefore building a same-size cube diagram from a given grid diagram cannot be done in general.  However, a grid diagram can always be used to build a larger-sized cube diagram.  In fact, given a grid diagram of size $n$ of a knot, a cube diagram of the same knot always exists of size at most
\begin{equation} n + 2(\mbox{\# of bad crossings}) + (\mbox{\# of twisted bends}). \label{gridtocube}
\end{equation}
The bad crossings mentioned above are those crossings in the $(y,z)-$ and $(z,x)-$projections that do not follow the convention given in the definition of a cube diagram.  For many grid diagrams, this leads to large cube diagrams 3 or 4 times the size of the original grid.  For example, a grid diagram for a 7 crossing knot of size 9 may produce a cube diagram of size 35 or more.  In this paper we find cube diagram representations for many 7 crossing knots of size 9.


\medskip

In order to derive the formula above and show how to improve upon it, we need a precise definition of cube diagrams, which is the content of the next section.

\bigskip
\section{Cube Diagrams of Knots}
\bigskip

Let $n$ be a positive integer and let $\Gamma$ be the cube $[0,n]\times [0,n]\times [0,n] \subset \mathbb{R}^3$ thought of as $3$-dimensional Cartesian grid, i.e., a grid with integer valued vertices.  A \textit{flat of $\Gamma$} is any cuboid (a right rectangular prism) with integer vertices in $\Gamma$ such that there are two orthogonal edges of length $n$ with the remaining orthogonal edge of length $1$.  A flat with an edge of length 1 that is parallel to the $x$-axis, $y$-axis, or $z$-axis is called an {\em $x$-flat}, {\em $y$-flat}, or {\em $z$-flat} respectively.

\medskip

The embedding of a link in the cube $\Gamma$ can be described as follows.  A marking is a labeled point in $\BR^3$ with half-integer coordinates.  Mark unit cubes of $\Gamma$ with either an $X$, $Y$, or $Z$ such that the following {\em marking conditions} hold:
\begin{itemize}
    \item each flat has exactly one $X$, one $Y$, and one $Z$ marking;\\

    \item the markings in each flat forms a right angle such that each ray is parallel to a coordinate axis;\\

    \item for each $x$-flat, $y$-flat, or $z$-flat, the marking that is the vertex of the right angle is an $X, Y,$ or $Z$ marking respectively.
\end{itemize}

\medskip

An oriented link can be embedded into $\Gamma$ by connecting pairs of markings with a line segment whenever two of their corresponding coordinates are the same.   Each line segment is oriented to go from an $X$ to a $Y$, from a $Y$ to a $Z$, or from a $Z$ to an $X$ (note that the cube itself is canonically oriented by the standard right hand orientation of $\BR^3$). The markings in each flat define two perpendicular segments of the link $L$ joined at a vertex, call the union of these segments a {\it cube bend}. If a cube bend is contained in an $x$-flat, we call it an {\it $x$-cube bend}. Similarly, define {\it $y$-cube bends} and {\it $z$-cube bends}.

\medskip

Arrange the markings in $\Gamma$ so that the following {\em crossing conditions} hold:
\begin{itemize}
\item At every intersection point of the $(x,y)$-projection, the segment parallel to the $x$-axis has smaller $z$-coordinate than the segment parallel to the $y$-axis.\\

\item At every intersection point of the $(y,z)$-projection, the segment parallel to the $y$-axis has smaller $x$-coordinate than the segment parallel to the $z$-axis.\\

\item At every intersection point of the $(z,x)$-projection, the segment parallel to the $z$-axis has smaller $y$-coordinate than the segment parallel to the $x$-axis.
\end{itemize}

If $\Gamma$ satisfies these conditions, then it is called a {\it cube diagram}.  We say that $\Gamma$ is a cube diagram representing the (oriented) link $L$.

\medskip

\begin{figure}[H]
\includegraphics[scale = .3]{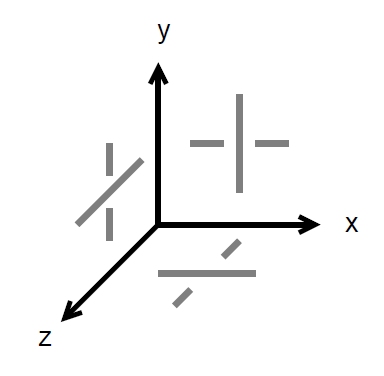}
\caption{Crossing conditions for the projections}
\label{crossings}
\end{figure}

The knot projections of a cube diagram to the three coordinate axis planes are grid diagrams.  What is not obvious, and the definition above specifies, is how the three projected grid diagrams are oriented with respect to the cube diagram.  In all three cases, the orientation is specified by the order of the axes using the standard orientation of $\BR^3$.  For example, in the $(x,y)$-projection, the $x$-axis specifies the `row' and the $y$-axis specifies the `column' of the grid diagram but in the $(y,z)$-projection the $y$-axis specifies the `row' and the $z$-axis specifies the `column' of that grid diagram (see Figure~\ref{crossings}).

These orientations matter:  if in the definition above, the first crossing condition was modified so that the segment parallel to the $x$-axis has {\em greater} $z$-coordinate than the segment parallel to the $y$-axis, then the set of knots types that are represented by modified `cube diagrams' of a given size is  different than the set of knot types using the actual definition.

\bigskip
\section{Building a cube diagram from a given grid diagram}
\bigskip

Section~2 described the problems with lifting a grid diagram to a cube diagram.  However, a grid diagram can always be used to create a cube diagram.  There are two issues to overcome: (1) changing the grid diagram so that it can be the knot-projection of a lattice knot such that the other two projections are well-defined knot projections and (2) fixing the crossings of that lattice knot in the other projections so that each crossing satisfies the crossing conditions.

\medskip

The first issue can be fixed by removing all of the twisted bends from a grid diagram.  A {\it bend} in a grid diagram for a knot $K$ is a pair of segments in $K$ that meet at a common $X$ or $O$ marking.  If a bend passes over some other segment of $K$ and passes under some other segment of $K$, then call it {\it twisted} (cf. \cite{Adam}).  There are two ways to partition a grid diagram of a knot into non-overlapping bends, depending on whether the two segments in each bend intersect in  an $X$ or $O$ marking.  A grid diagram of size $n$ can be used to construct a lattice knot embedded into an $n\times n\times n$ cube such that the other projections are knot projections if a partial order can be put on either of the two partitions of bends following the convention that if two bends cross the bend crossing over the other is greater (cf. \cite{Adam}).  If a grid diagram has a partition with no twisted bends, then a partial order always exists for that partition.  If a partition has one twisted bend, then stabilizing at the vertex of the twisted bend produces a new grid diagram with no twisted bends.  The new grid is one size larger, which explains the third term in the formula above when there are multiple twisted bends.

\medskip

\begin{figure}[H]
\includegraphics[scale = .3]{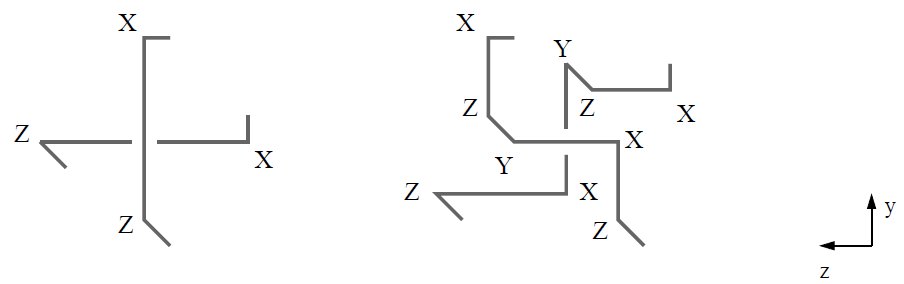}
\caption{Insertion of a rotated crossing}
\label{Methodology:Thor}
\end{figure}

If a grid diagram has no twisted bends, the bends can be stacked to form a lattice knot that projects to valid knot projections in all three planes.  While these knots do indeed satisfy the crossing conditions for the $(x,y)$-projection they may not satisfy the crossing conditions in one or both of the other projections.  As observed in \cite{Adam} the invalid crossings may be repaired by the insertion of a rotated crossing as shown in figure \ref{Methodology:Thor} (c.f. \cite{Adam}).

The rotated crossing corrects the crossing condition for that projection while leaving the crossing data of the other two projections unaffected.  This procedure has the effect of rotating the crossing so that the overcrossing is correct, but at the cost of increasing the size of the cube diagram by 2, which explains the second term in the formula above.

\bigskip
\section{Small cube diagrams for knots with small crossing knots}
\bigskip

The algorithm in Section 4 tends to produces large cube diagram representations for a given knot that are computationally intractable for computers to calculate invariants like knot Floer homology.  To find small cube diagrams for small crossing knots, we wrote a computer program that searches all size $\leq 9$ grid diagrams looking for grid diagrams that lift to cube diagrams.  If a grid diagram does, we check to see what knot type it is and record the cube diagram.  The idea is simple enough, but virtually all grid diagrams that lift to cube diagrams are unknots, and checking whether each grid diagram is possibly an unknot involves an $O(n^3)$ number of calculations per grid diagram of size $n$.  The total number of size $n$ grids (for knots and links) is

\smallskip

$$\frac{(n!)^2}{4}\left(1+2(1+n)\frac{\Gamma(1+n,-1)}{\euler \ \Gamma(2+n)}-6\frac{\Gamma(3,-1)}{\euler\ \Gamma(4)}\right),$$

\medskip

\noindent where $\Gamma(s)$ is the gamma function and $\Gamma(s,x)$ is the incomplete gamma function.  Clearly, checking whether each is an unknot is too time intensive for size 9 grids.  We describe next how we reduced the number of calculations to a routine that runs in days rather than years.

\medskip

Generating grid diagrams of size $n$ involves the choice of two size $n$ permutations $\sigma,\tau \in S_n$.  The coordinates of the $X$-markings of the grid diagram are given by $(i,\sigma(i))$ and the $O$-markings are given by $(i,\tau(i))$.  The program generates every possible grid diagram by an outer loop/inner loop structure.  The outer loop generates a new grid diagram by always starting with $\sigma=(1,2,3,\dots,n)$ ($X$'s along the diagonal) and cycling through the choices of $\tau$ that generate a grid diagram.  The inner loop then runs through all permutations of the columns using a routine that picks the `next largest' permutation in lexicographic order.  For example, if $n=3$, then the order of permutations from `smallest' to `largest' is $\sigma_1= (1,2,3), \sigma_2=(1,3,2), \sigma_3=(2,1,3), \sigma_4=(2,3,1), \sigma_5=(3,1,2),$ and $\sigma_6=(3,2,1)$.  By setting up the inner and outer loop in this way we reduce the number of times we need to check for the unknot significantly because the difference between the two diagrams from $\sigma_i$ and $\sigma_{i+1}$ is often either a (1) column commutation move or (2) a column cyclic permutation move.  Therefore if the grid diagram associated to $\sigma_i$ is the unknot and $\sigma_{i+1}$ is a commutation or cyclic permutation of $\sigma_i$, then we know that the grid diagram associated to $\sigma_{i+1}$ is also an unknot (no time intensive calculations are necessary).  Specifically, if the current grid diagram is merely a commutation of the previous diagram and the knot determinant of the previous diagram is 1, then the current diagram is also considered a potential unknot.  It is reasonable to throw out determinant 1 knots---according to \cite{knotinfo} there are just two nontrivial knots, $10_{124}$ and $12_{242}$, that have arc index $\leq 9$ and determinant 1.

\medskip

Once the program finds a grid diagram that potentially represents a nontrivial knot, the program looks for several $X-O$ configurations known not to lift to a cube diagram including eliminating links.  These configurations are relatively easy to check for and also significantly reduce the number of diagrams the program attempts to lift to cube diagrams.  Next, the knot determinant is computed.  As before, if the determinant is 1, the diagram is discarded as the unknot, $10_{124}$ knot, or $12_{242}$ knot.  The program then attempts to lift the grid to a valid cube by choosing a third permutation $\zeta \in S_n$ that determines the order in which the $z$-cube bends will be stacked in order from smallest $z$-coordinate to greatest.  For each permutation the program checks to ensure that the stack is compatible with the crossing conditions determined by the grid.  If a valid stack permutation is found, the crossing conditions are checked in the $(y,z)-$ and $(z,x)$-projections.  If the crossing conditions are satisfied, the tuple of permutations $(\sigma,\tau,\zeta)$ is a cube diagram and the program computes a variant of the Jones polynomial to determine exactly what knot the cube diagram represents (c.f. \cite{Zulli}).  This final and most time intensive calculation involves $O(c^22^c)$ number of calculations where $c$ is the crossing number of the diagram.  Fortunately, this calculation is not necessary very often.

Using the routine sketched out above and beginning with grid diagrams of size $\leq 9$ we obtained most of the list in Appendix \ref{appendix:examples}.  For the 6 crossing knots, however, a slightly different strategy was used.  Beginning with a valid grid diagram for each knot type the cube stacking algorithm described above was used to produce lattice knots that projected to valid knot projections in all three planes.  While these knots did indeed satisfy the crossing conditions for the $(x,y)$-projection they all had invalid crossings in at least one of the other two projections.  As observed above, the invalid crossings may be repaired by rotating the crossing as shown in Figure \ref{Methodology:Thor}.  It is not known if these examples are the smallest.

\appendix
\section{Knot Examples}
\label{appendix:examples}
\medskip
The following table lists examples of cube knots up to 7 crossings as well as several 8, 9, 10 and 12 crossing knots.  Each picture displays the cube knot from the point of view of the $(x,y)$-projection.  The code presented to the right of each diagram is designed to work with a Mathematica notebook found at \cite{BLcode}.  Since projections of cube diagrams are grid diagrams, the arc index $\alpha(K)$ of a knot gives a lower bound for the size of a cube diagram.  It has been shown in \cite{arcindex} that for alternating knots $\alpha(K) = c(K) + 2$ where $c(K)$ is the crossing number of the knot.  Therefore in searching for cube diagrams up to size 9 the program could only be expected to find alternating knots up to 7 crossings (one exception below, $8_{15}$, was found by running a partial search for size 10 diagrams).  However, for non-alternating knots the arc index may be much smaller.  This fact helps to explain why some of the 9, 10 and 12 crossing knots show up in the table below.  The examples that were found were of knots with relatively low arc index.  For example $12_{591}$ listed below is a non-alternating knot that has arc index equal to 9 (see \cite{knotinfo}).

It is interesting to note that for the 6-crossing knots, $7_6$, and $7_7$ the cube diagrams were obtained using the second method described above.  This method was necessary because the program found no valid cube diagrams of size $\leq 9$.

\medskip

\begin{center}

\begin{tabular}{| c | l |}
\hline
Diagram & Mathematica Code \\
\hline
\multirow{10}{*}{\includegraphics[scale = .3]{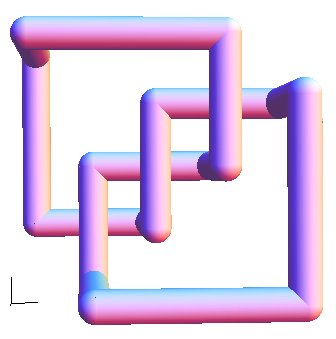}} & \\
& \\
& \\
& K$3_1$ = \{X[\{1, 5, 4\}, \{4, 3, 2\}, \{5, 4, 3\}, \{2, 1, 5\}, \{3, 2, 1\}], \\
& Y[\{1, 5, 1\}, \{2, 1, 2\}, \{3, 2, 3\}, \{4, 3, 4\}, \{5, 4, 5\}], \\
& Z[\{1, 2, 1\}, \{2, 3, 2\}, \{3, 4, 3\}, \{4, 5, 4\}, \{5, 1, 5\}]\} \\
& \\
& \\
& \\
$3_1$ & \\
\hline
\multirow{11}{*}{\includegraphics[scale = .3]{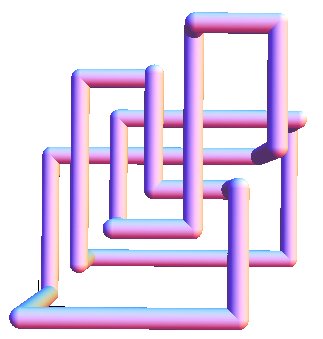}} & \\
& \\
& K$4_1$ = \{X[\{1, 1, 2\}, \{2, 2, 4\}, \{3, 3, 3\}, \{4, 7, 5\}, \{5, 8, 6\}, \\
& \{6, 4, 8\}, \{7, 5, 7\}, \{8, 6, 1\}],\\
& Y[\{1, 1, 8\}, \{2, 2, 1\}, \{3, 3, 6\}, \{4, 7, 4\}, \{5, 8, 7\}, \\
& \{6, 4, 5\}, \{7, 5, 2\}, \{8, 6, 3\}],\\
& Z[\{1, 5, 2\}, \{2, 7, 4\}, \{3, 6, 3\}, \{4, 4, 5\}, \{5, 3, 6\}, \\
& \{6, 1, 8\}, \{7, 8, 7\}, \{8, 2, 1\}]\}\\
& \\
& \\
$4_1$ & \\
\hline
\end{tabular}

\begin{tabular}{| c | l |}
\hline
Diagram & Mathematica Code \\
\hline
\multirow{11}{*}{\includegraphics[scale = .3]{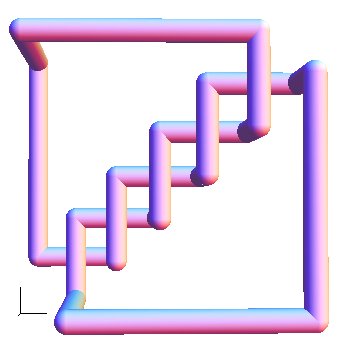}} & \\
& \\
& K$5_1$ = \{X[\{1, 7, 1\}, \{2, 1, 2\}, \{3, 2, 3\}, \{4, 3, 4\}, \{5, 4, 5\}, \\
& \{6, 5, 6\}, \{7, 6, 7\}],\\
& Y[\{1, 7, 6\}, \{2, 1, 7\}, \{3, 2, 1\}, \{4, 3, 2\}, \{5, 4, 3\}, \\
& \{6, 5, 4\}, \{7, 6, 5\}],\\
& Z[\{1, 2, 1\}, \{2, 3, 2\}, \{3, 4, 3\}, \{4, 5, 4\}, \{5, 6, 5\}, \\
& \{6, 7, 6\}, \{7, 1, 7\}]\}\\
& \\
& \\
$5_1$ & \\
\hline
\multirow{11}{*}{\includegraphics[scale = .3]{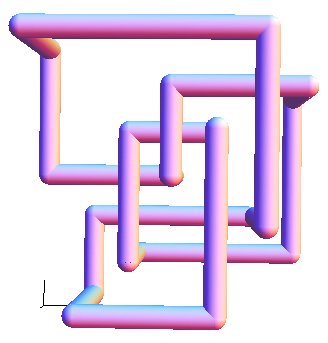}} & \\
& \\
& K$5_2$ = \{X[\{2, 1, 6\}, \{1, 7, 7\}, \{3, 2, 1\}, \{4, 4, 3\}, \{5, 5, 4\}, \\
& \{6, 3, 2\}, \{7, 6, 5\}],\\
& Y[\{2, 1, 2\}, \{1, 7, 3\}, \{3, 2, 4\}, \{4, 4, 5\}, \{5, 5, 6\}, \\
& \{6, 3, 7\}, \{7, 6, 1\}],\\
& Z[\{2, 3, 2\}, \{1, 4, 3\}, \{3, 5, 4\}, \{4, 6, 5\}, \{5, 1, 6\}, \\
& \{6, 7, 7\}, \{7, 2, 1\}]\}\\
& \\
& \\
$5_2$ & \\
\hline
\multirow{11}{*}{\includegraphics[scale = .3]{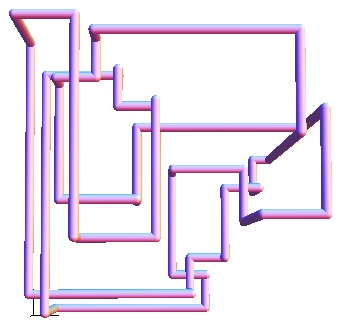}} & \\
& K$6_1$ = \{X[\{1, 16, 7\}, \{2, 1, 8\}, \{3, 13, 11\}, \{4, 5, 14\}, \{5, 15, 13\}, \\
& \{6, 14, 10\}, \{7, 7, 9\}, \{8, 12, 12\}, \{9, 9, 3\}, \{10, 2, 4\}, \\
& \{11, 3, 2\}, \{12, 4, 5\}, \{13, 6, 6\}, \{14, 8, 1\}, \{15, 10, 15\}, \\
& \{16, 11, 16\}],\\
& Y[\{1, 16, 14\}, \{2, 1, 2\}, \{3, 13, 13\}, \{4, 5, 12\}, \{5, 15, 15\}, \\
& \{6, 14, 8\}, \{7, 7, 11\}, \{8, 12, 10\}, \{9, 9, 6\}, \{10, 2, 7\}, \\
& \{11, 3, 3\}, \{12, 4, 4\}, \{13, 6, 16\}, \{14, 8, 5\}, \{15, 10, 1\}, \\
& \{16, 11, 9\}],\\
& Z[\{1, 2, 7\}, \{2, 14, 8\}, \{3, 7, 11\}, \{4, 16, 14\}, \{5, 13, 13\}, \\
& \{6, 12, 10\}, \{7, 11, 9\}, \{8, 5, 12\}, \{9, 3, 3\}, \{10, 4, 4\}, \\
& \{11, 1, 2\}, \{12, 8, 5\}, \{13, 9, 6\}, \{14, 10, 1\}, \{15, 15, 15\}, \\
& \{16, 6, 16\}]\}\\
$6_1$ & \\
\hline
\multirow{11}{*}{\includegraphics[scale = .3]{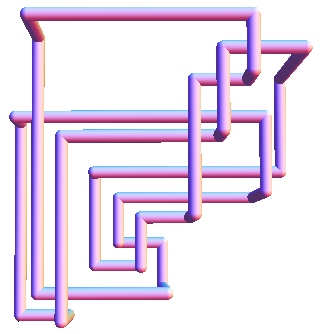}} & \\
& K$6_2$ = \{X[\{1, 9, 5\}, \{2, 12, 8\}, \{3, 1, 9\}, \{4, 7, 3\}, \{5, 4, 2\}, \\
& \{6, 3, 6\}, \{7, 2, 1\}, \{8, 5, 10\}, \{9, 8, 11\}, \{10, 10, 12\}, \\
& \{11, 6, 7\}, \{12, 11, 4\}],\\
& Y[\{1, 9, 7\}, \{2, 12, 12\}, \{3, 1, 5\}, \{4, 7, 4\}, \{5, 4, 1\}, \\
& \{6, 3, 3\}, \{7, 2, 8\}, \{8, 5, 6\}, \{9, 8, 9\}, \{10, 10, 10\}, \\
& \{11, 6, 2\}, \{12, 11, 11\}],\\
& Z[\{1, 1, 5\}, \{2, 2, 8\}, \{3, 8, 9\}, \{4, 3, 3\}, \{5, 6, 2\}, \\
& \{6, 5, 6\}, \{7, 4, 1\}, \{8, 10, 10\}, \{9, 11, 11\}, \{10, 12, 12\}, \\
& \{11, 9, 7\}, \{12, 7, 4\}]\}\\
$6_2$ & \\
\hline
\end{tabular}

\begin{tabular}{| c | l |}
\hline
Diagram & Mathematica Code \\
\hline
\multirow{11}{*}{\includegraphics[scale = .3]{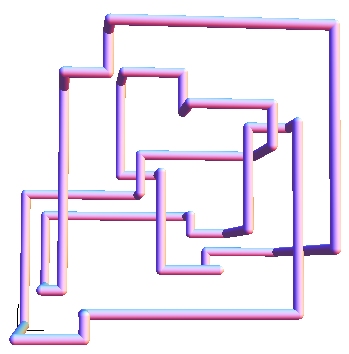}} & \\
& K$6_3$ = \{X[\{1, 1, 4\}, \{2, 3, 1\}, \{3, 13, 9\}, \{4, 2, 14\}, \{5, 15, 13\}, \\
& \{6, 14, 5\}, \{7, 8, 6\}, \{8, 9, 7\}, \{9, 12, 8\}, \{10, 7, 3\}, \\
& \{11, 4, 2\}, \{12, 6, 10\}, \{13, 10, 11\}, \{14, 11, 12\}, \{15, 5, 15\}],\\
& Y[\{1, 1, 14\}, \{2, 3, 9\}, \{3, 13, 13\}, \{4, 2, 12\}, \{5, 15, 15\}, \\
& \{6, 14, 8\}, \{7, 8, 4\}, \{8, 9, 5\}, \{9, 12, 11\}, \{10, 7, 1\}, \\
& \{11, 4, 7\}, \{12, 6, 3\}, \{13, 10, 6\}, \{14, 11, 10\}, \{15, 5, 2\}],\\
& Z[\{1, 8, 4\}, \{2, 7, 1\}, \{3, 3, 9\}, \{4, 1, 14\}, \{5, 13, 13\}, \\
& \{6, 9, 5\}, \{7, 10, 6\}, \{8, 4, 7\}, \{9, 14, 8\}, \{10, 6, 3\}, \\
& \{11, 5, 2\}, \{12, 11, 10\}, \{13, 12, 11\}, \{14, 2, 12\}, \{15, 15, 15\}]\}\\
$6_3$ & \\
\hline
\multirow{11}{*}{\includegraphics[scale = .3]{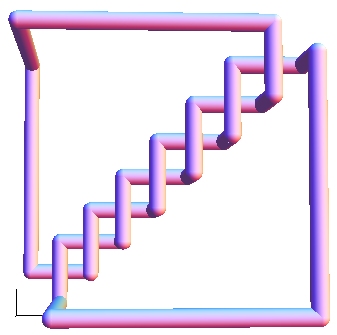}} & \\
& \\
& K$7_1$ = \{X[\{1, 9, 1\}, \{2, 1, 2\}, \{3, 2, 3\}, \{4, 3, 4\}, \{5, 4, 5\}, \\
& \{6, 5, 6\}, \{7, 6, 7\}, \{8, 7, 8\}, \{9, 8, 9\}],\\
& Y[\{1, 9, 8\}, \{2, 1, 9\}, \{3, 2, 1\}, \{4, 3, 2\}, \{5, 4, 3\}, \\
& \{6, 5, 4\}, \{7, 6, 5\}, \{8, 7, 6\}, \{9, 8, 7\}],\\
& Z[\{1, 2, 1\}, \{2, 3, 2\}, \{3, 4, 3\}, \{4, 5, 4\}, \{5, 6, 5\}, \\
& \{6, 7, 6\}, \{7, 8, 7\}, \{8, 9, 8\}, \{9, 1, 9\}]\}\\
& \\
& \\
$7_1$ & \\
\hline
\multirow{11}{*}{\includegraphics[scale = .3]{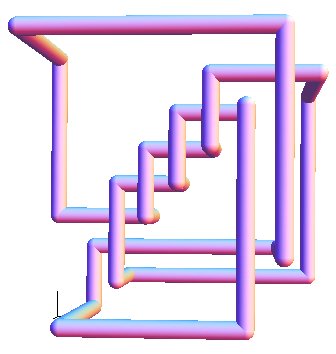}} & \\
& \\
& K$7_2$ = \{X[\{2, 1, 2\}, \{1, 9, 3\}, \{3, 2, 4\}, \{4, 4, 5\}, \{5, 5, 6\}, \\
& \{6, 6, 7\}, \{7, 7, 8\}, \{8, 3, 9\}, \{9, 8, 1\}],\\
& Y[\{2, 1, 8\}, \{1, 9, 9\}, \{3, 2, 1\}, \{4, 4, 3\}, \{5, 5, 4\}, \\
& \{6, 6, 5\}, \{7, 7, 6\}, \{8, 3, 2\}, \{9, 8, 7\}],\\
& Z[\{2, 3, 2\}, \{1, 4, 3\}, \{3, 5, 4\}, \{4, 6, 5\}, \{5, 7, 6\}, \\
& \{6, 8, 7\}, \{7, 1, 8\}, \{8, 9, 9\}, \{9, 2, 1\}]\}\\
& \\
& \\
$7_2$ & \\
\hline
\multirow{11}{*}{\includegraphics[scale = .3]{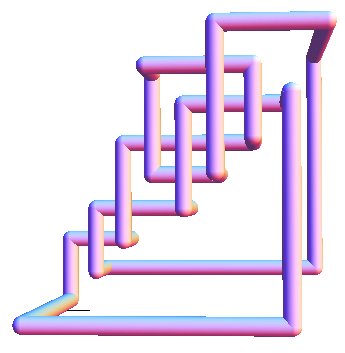}} & \\
& \\
& K$7_3$ = \{X[\{1, 1, 2\}, \{2, 2, 3\}, \{3, 3, 4\}, \{4, 8, 5\}, \{5, 4, 6\}, \\
& \{6, 5, 8\}, \{7, 6, 7\}, \{8, 7, 9\}, \{9, 9, 1\}],\\
& Y[\{1, 1, 9\}, \{2, 2, 1\}, \{3, 3, 2\}, \{4, 8, 7\}, \{5, 4, 3\}, \\
& \{6, 5, 5\}, \{7, 6, 4\}, \{8, 7, 6\}, \{9, 9, 8\}],\\
& Z[\{1, 3, 2\}, \{2, 4, 3\}, \{3, 6, 4\}, \{4, 5, 5\}, \{5, 7, 6\}, \\
& \{6, 9, 8\}, \{7, 8, 7\}, \{8, 1, 9\}, \{9, 2, 1\}]\}\\
& \\
& \\
$7_3$ & \\
\hline
\end{tabular}

\begin{tabular}{| c | l |}
\hline
Diagram & Mathematica Code \\
\hline
\multirow{11}{*}{\includegraphics[scale = .3]{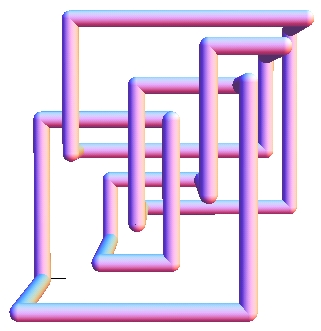}} & \\
& \\
& K$7_4$ = \{X[\{1, 1, 4\}, \{2, 5, 5\}, \{3, 2, 2\}, \{4, 3, 6\}, \{5, 6, 7\}, \\
& \{6, 4, 8\}, \{7, 7, 9\}, \{8, 8, 3\}, \{9, 9, 1\}],\\
& Y[\{1, 1, 9\}, \{2, 5, 3\}, \{3, 2, 7\}, \{4, 3, 1\}, \{5, 6, 4\}, \\
& \{6, 4, 2\}, \{7, 7, 6\}, \{8, 8, 8\}, \{9, 9, 5\}],\\
& Z[\{1, 6, 4\}, \{2, 9, 5\}, \{3, 4, 2\}, \{4, 7, 6\}, \{5, 2, 7\}, \\
& \{6, 8, 8\}, \{7, 1, 9\}, \{8, 5, 3\}, \{9, 3, 1\}]\}\\
& \\
& \\
$7_4$ & \\
\hline
\multirow{11}{*}{\includegraphics[scale = .3]{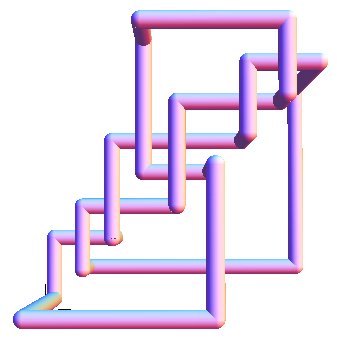}} & \\
& \\
& K$7_5$ = \{X[\{1, 1, 2\}, \{2, 2, 3\}, \{3, 3, 4\}, \{4, 9, 5\}, \{5, 4, 6\}, \\
& \{6, 5, 9\}, \{7, 6, 7\}, \{8, 7, 8\}, \{9, 8, 1\}],\\
& Y[\{1, 1, 9\}, \{2, 2, 1\}, \{3, 3, 2\}, \{4, 9, 8\}, \{5, 4, 3\}, \\
& \{6, 5, 5\}, \{7, 6, 4\}, \{8, 7, 6\}, \{9, 8, 7\}],\\
& Z[\{1, 3, 2\}, \{2, 4, 3\}, \{3, 6, 4\}, \{4, 5, 5\}, \{5, 7, 6\}, \\
& \{6, 1, 9\}, \{7, 8, 7\}, \{8, 9, 8\}, \{9, 2, 1\}]\}\\
& \\
& \\
$7_5$ & \\
\hline
\multirow{11}{*}{\includegraphics[scale = .3]{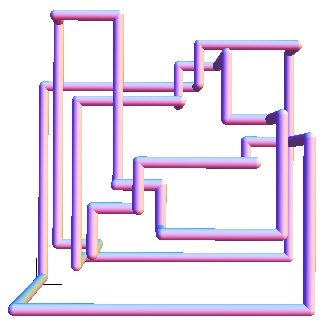}} & \\
& K$7_6$ = \{X[\{1, 1, 2\}, \{2, 14, 5\}, \{3, 2, 6\}, \{4, 3, 11\}, \{5, 6, 7\}, \\
& \{6, 5, 12\}, \{7, 4, 10\}, \{8, 10, 3\}, \{9, 11, 4\}, \{10, 12, 9\}, \\
& \{11, 7, 8\}, \{12, 9, 13\}, \{13, 8, 14\}, \{14, 13, 1\}],\\
& Y[\{1, 1, 14\}, \{2, 14, 7\}, \{3, 2, 1\}, \{4, 3, 5\}, \{5, 6, 10\}, \\
& \{6, 5, 11\}, \{7, 4, 13\}, \{8, 10, 6\}, \{9, 11, 2\}, \{10, 12, 3\}, \\
& \{11, 7, 12\}, \{12, 9, 9\}, \{13, 8, 8\}, \{14, 13, 4\}],\\
& Z[\{1, 11, 2\}, \{2, 3, 5\}, \{3, 10, 6\}, \{4, 5, 11\}, \{5, 14, 7\}, \\
& \{6, 7, 12\}, \{7, 6, 10\}, \{8, 12, 3\}, \{9, 13, 4\}, \{10, 9, 9\}, \\
& \{11, 8, 8\}, \{12, 4, 13\}, \{13, 1, 14\}, \{14, 2, 1\}]\}\\
$7_6$ & \\
\hline
\multirow{13}{*}{\includegraphics[scale = .3]{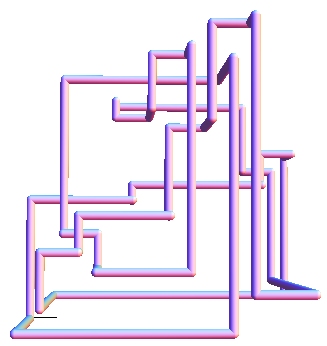}} & \\
& K$7_7$ = \{X[\{1, 1, 7\}, \{2, 2, 13\}, \{3, 6, 8\}, \{4, 5, 12\}, \{5, 4, 9\}, \\
& \{6, 14, 6\}, \{7, 8, 5\}, \{8, 13, 11\}, \{9, 7, 10\}, \{10, 16, 14\}, \\
& \{11, 12, 15\}, \{12, 15, 16\}, \{13, 17, 17\}, \{14, 10, 3\}, \{15, 9, 4\}, \\
& \{16, 3, 2\}, \{17, 11, 1\}],\\
& Y[\{1, 1, 16\}, \{2, 2, 1\}, \{3, 6, 9\}, \{4, 5, 13\}, \{5, 4, 14\}, \\
& \{6, 14, 3\}, \{7, 8, 7\}, \{8, 13, 6\}, \{9, 7, 12\}, \{10, 16, 11\}, \\
& \{11, 12, 10\}, \{12, 15, 8\}, \{13, 17, 15\}, \{14, 10, 2\}, \{15, 9, 5\}, \\
& \{16, 3, 17\}, \{17, 11, 4\}],\\
& Z[\{1, 8, 7\}, \{2, 5, 13\}, \{3, 15, 8\}, \{4, 7, 12\}, \{5, 6, 9\}, \\
& \{6, 13, 6\}, \{7, 9, 5\}, \{8, 16, 11\}, \{9, 12, 10\}, \{10, 4, 14\}, \\
& \{11, 17, 15\}, \{12, 1, 16\}, \{13, 3, 17\}, \{14, 14, 3\}, \{15, 11, 4\}, \\
& \{16, 10, 2\}, \{17, 2, 1\}]\}\\
$7_7$ & \\
\hline
\end{tabular}

\begin{tabular}{| c | l |}
\hline
Diagram & Mathematica Code \\
\hline
\multirow{11}{*}{\includegraphics[scale = .3]{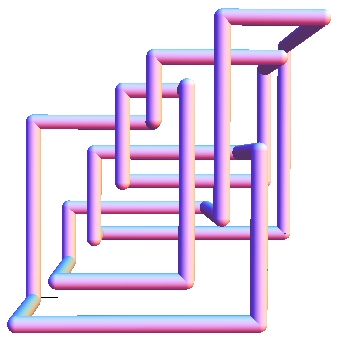}} & \\
& \\
& K$8_{15}$ = \{X[\{1, 1, 5\}, \{2, 2, 2\}, \{3, 3, 4\}, \{4, 5, 6\}, \{5, 7, 7\}, \\
& \{6, 8, 8\}, \{7, 4, 9\}, \{8, 6, 10\}, \{9, 9, 3\}, \{10, 10, 1\}],\\
& Y[\{1, 1, 10\}, \{2, 2, 8\}, \{3, 3, 1\}, \{4, 5, 3\}, \{5, 7, 5\}, \\
& \{6, 8, 6\}, \{7, 4, 2\}, \{8, 6, 4\}, \{9, 9, 7\}, \{10, 10, 9\}],\\
& Z[\{1, 7, 5\}, \{2, 4, 2\}, \{3, 6, 4\}, \{4, 8, 6\}, \{5, 9, 7\}, \\
& \{6, 2, 8\}, \{7, 10, 9\}, \{8, 1, 10\}, \{9, 5, 3\}, \{10, 3, 1\}]\}\\
& \\
& \\
$8_{15}$ & \\
\hline
\multirow{11}{*}{\includegraphics[scale = .3]{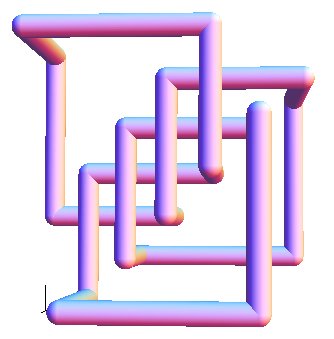}} & \\
& \\
& K$8_{19}$ = \{X[\{1, 7, 2\}, \{2, 1, 3\}, \{3, 2, 4\}, \{4, 3, 5\}, \{5, 4, 6\}, \\
& \{6, 5, 7\}, \{7, 6, 1\}],\\
& Y[\{1, 7, 6\}, \{2, 1, 7\}, \{3, 2, 1\}, \{4, 3, 2\}, \{5, 4, 3\}, \\
& \{6, 5, 4\}, \{7, 6, 5\}],\\
& Z[\{1, 3, 2\}, \{2, 4, 3\}, \{3, 5, 4\}, \{4, 6, 5\}, \{5, 7, 6\}, \\
& \{6, 1, 7\}, \{7, 2, 1\}]\}\\
& \\
& \\
$8_{19}$ & \\
\hline
\multirow{11}{*}{\includegraphics[scale = .3]{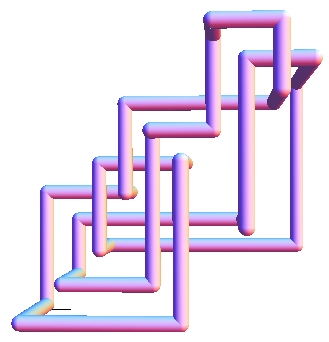}} & \\
& \\
& K$8_{21}$ = \{X[\{1, 1, 3\}, \{2, 2, 2\}, \{3, 3, 4\}, \{4, 5, 5\}, \{5, 7, 6\}, \\
& \{6, 6, 8\}, \{7, 10, 7\}, \{8, 4, 9\}, \{9, 8, 10\}, \{10, 9, 1\}],\\
& Y[\{1, 1, 8\}, \{2, 2, 6\}, \{3, 3, 1\}, \{4, 5, 3\}, \{5, 7, 7\}, \\
& \{6, 6, 4\}, \{7, 10, 10\}, \{8, 4, 2\}, \{9, 8, 5\}, \{10, 9, 9\}],\\
& Z[\{1, 5, 3\}, \{2, 4, 2\}, \{3, 6, 4\}, \{4, 8, 5\}, \{5, 2, 6\}, \\
& \{6, 1, 8\}, \{7, 7, 7\}, \{8, 9, 9\}, \{9, 10, 10\}, \{10, 3, 1\}]\}\\
& \\
& \\
$8_{21}$ & \\
\hline
\multirow{11}{*}{\includegraphics[scale = .3]{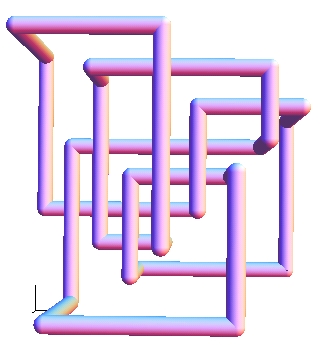}} & \\
& \\
& K$9_{49}$ = \{X[\{1, 9, 2\}, \{2, 1, 3\}, \{3, 8, 4\}, \{4, 2, 5\}, \{5, 3, 8\}, \\
& \{6, 4, 6\}, \{7, 5, 9\}, \{8, 6, 7\}, \{9, 7, 1\}],\\
& Y[\{1, 9, 8\}, \{2, 1, 9\}, \{3, 8, 7\}, \{4, 2, 1\}, \{5, 3, 4\}, \\
& \{6, 4, 2\}, \{7, 5, 5\}, \{8, 6, 3\}, \{9, 7, 6\}],\\
& Z[\{1, 4, 2\}, \{2, 6, 3\}, \{3, 3, 4\}, \{4, 5, 5\}, \{5, 9, 8\}, \\
& \{6, 7, 6\}, \{7, 1, 9\}, \{8, 8, 7\}, \{9, 2, 1\}]\}\\
& \\
& \\
$9_{49}$ & \\
\hline
\end{tabular}

\begin{tabular}{| c | l |}
\hline
Diagram & Mathematica Code \\
\hline
\multirow{11}{*}{\includegraphics[scale = .3]{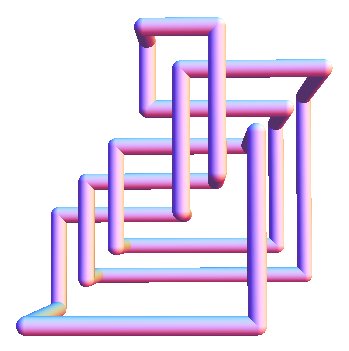}} & \\
& \\
& K$10_{124}$ = \{X[\{1, 1, 3\}, \{2, 2, 4\}, \{3, 3, 5\}, \{4, 9, 6\}, \{5, 4, 7\}, \\
& \{6, 5, 8\}, \{7, 6, 9\}, \{8, 7, 2\}, \{9, 8, 1\}],\\
& Y[\{1, 1, 9\}, \{2, 2, 1\}, \{3, 3, 2\}, \{4, 9, 8\}, \{5, 4, 3\}, \\
& \{6, 5, 4\}, \{7, 6, 5\}, \{8, 7, 6\}, \{9, 8, 7\}],\\
& Z[\{1, 4, 3\}, \{2, 5, 4\}, \{3, 6, 5\}, \{4, 7, 6\}, \{5, 8, 7\}, \\
& \{6, 9, 8\}, \{7, 1, 9\}, \{8, 3, 2\}, \{9, 2, 1\}]\}\\
& \\
& \\
$10_{124}$ & \\
\hline
\multirow{11}{*}{\includegraphics[scale = .3]{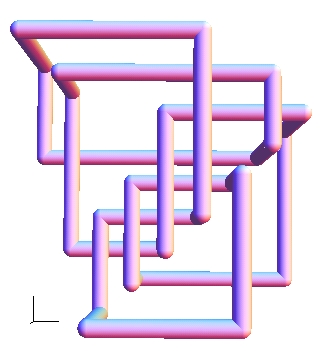}} & \\
& \\
& K$10_{128}$ = \{X[\{1, 9, 2\}, \{2, 8, 3\}, \{3, 1, 4\}, \{4, 2, 5\}, \{5, 3, 6\}, \\
& \{6, 4, 8\}, \{7, 5, 9\}, \{8, 6, 7\}, \{9, 7, 1\}],\\
& Y[\{1, 9, 8\}, \{2, 8, 7\}, \{3, 1, 9\}, \{4, 2, 1\}, \{5, 3, 3\}, \\
& \{6, 4, 4\}, \{7, 5, 5\}, \{8, 6, 2\}, \{9, 7, 6\}],\\
& Z[\{1, 6, 2\}, \{2, 3, 3\}, \{3, 4, 4\}, \{4, 5, 5\}, \{5, 7, 6\}, \\
& \{6, 9, 8\}, \{7, 1, 9\}, \{8, 8, 7\}, \{9, 2, 1\}]\}\\
& \\
& \\
$10_{128}$ & \\
\hline
\multirow{11}{*}{\includegraphics[scale = .3]{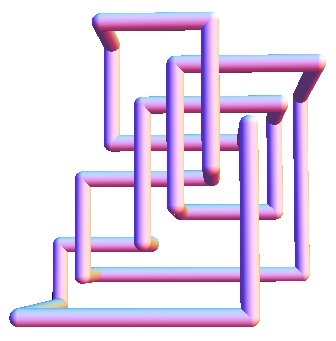}} & \\
& \\
& K$10_{139}$ = \{X[\{1, 1, 2\}, \{2, 2, 4\}, \{3, 9, 5\}, \{4, 3, 6\}, \{5, 4, 7\}, \\
& \{6, 5, 8\}, \{7, 6, 9\}, \{8, 7, 3\}, \{9, 8, 1\}],\\
& Y[\{1, 1, 9\}, \{2, 2, 1\}, \{3, 9, 8\}, \{4, 3, 2\}, \{5, 4, 3\}, \\
& \{6, 5, 4\}, \{7, 6, 5\}, \{8, 7, 6\}, \{9, 8, 7\}],\\
& Z[\{1, 3, 2\}, \{2, 5, 4\}, \{3, 6, 5\}, \{4, 7, 6\}, \{5, 8, 7\}, \\
& \{6, 9, 8\}, \{7, 1, 9\}, \{8, 4, 3\}, \{9, 2, 1\}]\}\\
& \\
& \\
$10_{139}$ & \\
\hline
\multirow{11}{*}{\includegraphics[scale = .3]{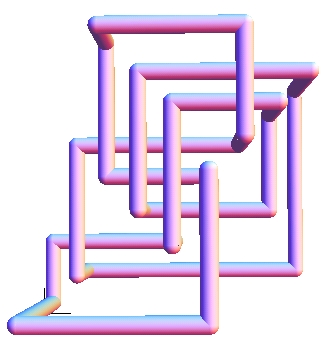}} & \\
& \\
& K$10_{145}$ = \{X[\{1, 1, 2\}, \{2, 2, 4\}, \{3, 9, 5\}, \{4, 4, 6\}, \{5, 3, 7\}, \\
& \{6, 5, 9\}, \{7, 6, 8\}, \{8, 7, 3\}, \{9, 8, 1\}],\\
& Y[\{1, 1, 9\}, \{2, 2, 1\}, \{3, 9, 8\}, \{4, 4, 3\}, \{5, 3, 2\}, \\
& \{6, 5, 5\}, \{7, 6, 4\}, \{8, 7, 7\}, \{9, 8, 6\}],\\
& Z[\{1, 3, 2\}, \{2, 6, 4\}, \{3, 5, 5\}, \{4, 8, 6\}, \{5, 7, 7\}, \\
& \{6, 1, 9\}, \{7, 9, 8\}, \{8, 4, 3\}, \{9, 2, 1\}]\}\\
& \\
& \\
$10_{145}$ & \\
\hline
\end{tabular}

\begin{tabular}{| c | l |}
\hline
Diagram & Mathematica Code \\
\hline
\multirow{11}{*}{\includegraphics[scale = .3]{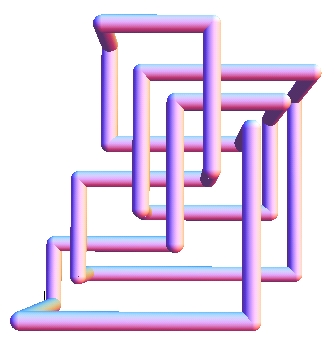}} & \\
& \\
& K$10_{161}$ = \{X[\{1, 1, 2\}, \{2, 2, 4\}, \{3, 9, 5\}, \{4, 4, 6\}, \{5, 3, 7\}, \\
& \{6, 5, 8\}, \{7, 6, 9\}, \{8, 7, 3\}, \{9, 8, 1\}],\\
& Y[\{1, 1, 9\}, \{2, 2, 1\}, \{3, 9, 8\}, \{4, 4, 3\}, \{5, 3, 2\}, \\
& \{6, 5, 4\}, \{7, 6, 5\}, \{8, 7, 7\}, \{9, 8, 6\}],\\
& Z[\{1, 3, 2\}, \{2, 5, 4\}, \{3, 6, 5\}, \{4, 8, 6\}, \{5, 7, 7\}, \\
& \{6, 9, 8\}, \{7, 1, 9\}, \{8, 4, 3\}, \{9, 2, 1\}]\}\\
& \\
& \\
$10_{161}$ & \\
\hline
\multirow{11}{*}{\includegraphics[scale = .3]{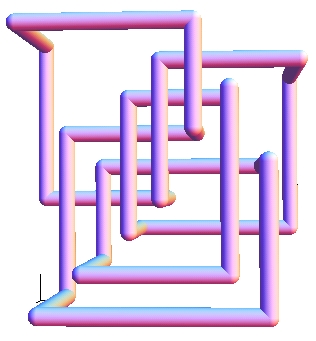}} & \\
& \\
& K$12_{591}$ = \{X[\{1, 9, 2\}, \{2, 1, 4\}, \{3, 2, 3\}, \{4, 3, 5\}, \{5, 4, 6\}, \\
& \{6, 6, 7\}, \{7, 7, 8\}, \{8, 5, 9\}, \{9, 8, 1\}],\\
& Y[\{1, 9, 7\}, \{2, 1, 9\}, \{3, 2, 8\}, \{4, 3, 1\}, \{5, 4, 2\}, \\
& \{6, 6, 4\}, \{7, 7, 5\}, \{8, 5, 3\}, \{9, 8, 6\}],\\
& Z[\{1, 4, 2\}, \{2, 6, 4\}, \{3, 5, 3\}, \{4, 7, 5\}, \{5, 8, 6\}, \\
& \{6, 9, 7\}, \{7, 2, 8\}, \{8, 1, 9\}, \{9, 3, 1\}]\}\\
& \\
& \\
$12_{591}$ & \\
\hline
\end{tabular}

\end{center}

\section{Link Examples}

\begin{center}

\begin{tabular}{| c | l |}
\hline
Diagram & Mathematica Code \\
\hline
& \\
\multirow{4}{*}{\includegraphics[scale = .2]{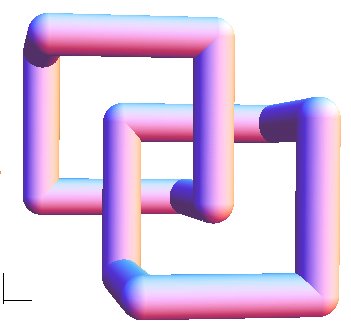}} & \\
& HL = \{X[\{1, 4, 1\}, \{2, 1, 2\}, \{3, 2, 3\}, \{4, 3, 4\}],\\
& Y[\{1, 4, 3\}, \{2, 1, 4\}, \{3, 2, 1\}, \{4, 3, 2\}],\\
& Z[\{1, 2, 1\}, \{2, 3, 2\}, \{3, 4, 3\}, \{4, 1, 4\}]\}\\
& \\
& \\
$Hopf Link$ & \\
\hline
\multirow{14}{*}{\includegraphics[scale = .3]{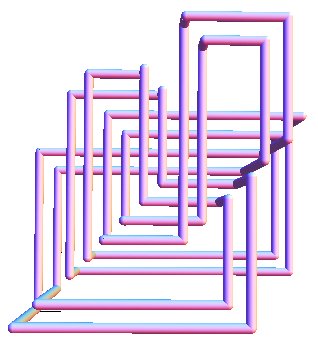}} & \\
& $D_{0}(4_1)$ = \{X[\{1, 1, 16\}, \{2, 2, 15\}, \{3, 3, 1\}, \{4, 4, 2\}, \{5, 5, 12\}, \\
& \{6, 6, 11\}, \{7, 14, 8\}, \{8, 13, 7\}, \{10, 15, 13\}, \{9, 16, 14\}, \\
& \{11, 7, 10\}, \{12, 8, 9\}, \{14, 10, 4\}, \{13, 9, 3\}, \{15, 11, 5\}, \\
& \{16, 12, 6\}],\\
& Y[\{1, 1, 4\}, \{2, 2, 3\}, \{3, 3, 7\}, \{4, 4, 8\}, \{5, 5, 6\}, \\
& \{6, 6, 5\}, \{7, 14, 10\}, \{8, 13, 9\}, \{10, 15, 11\}, \{9, 16, 12\}, \\
& \{11, 7, 15\}, \{12, 8, 16\}, \{14, 10, 14\}, \{13, 9, 13\}, \{15, 11, 2\}, \\
& \{16, 12, 1\}],\\
& Z[\{1, 10, 4\}, \{2, 9, 3\}, \{3, 13, 7\}, \{4, 14, 8\}, \{5, 12, 6\}, \\
& \{6, 11, 5\}, \{7, 7, 10\}, \{8, 8, 9\}, \{10, 6, 11\}, \{9, 5, 12\}, \\
& \{11, 2, 15\}, \{12, 1, 16\}, \{14, 16, 14\}, \{13, 15, 13\}, \{15, 4, 2\}, \\
& \{16, 3, 1\}]\}\\
$D_{0}(4_1)$ & \\
\hline
\multirow{11}{*}{\includegraphics[scale = .3]{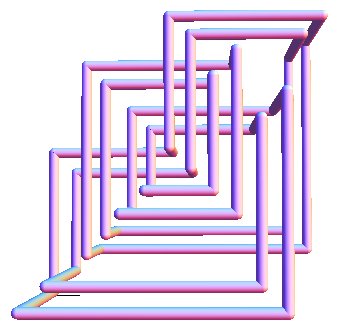}} & \\
& $D_{-4}(5_2)$ = \{X[\{1, 1, 4\}, \{2, 2, 3\}, \{3, 3, 7\}, \{4, 4, 8\}, \{5, 5, 5\}, \\
& \{6, 6, 6\}, \{7, 8, 9\}, \{8, 7, 10\}, \{9, 11, 11\}, \{10, 12, 12\}, \\
& \{11, 9, 13\}, \{12, 10, 14\}, \{13, 13, 2\}, \{14, 14, 1\}],\\
& Y[\{1, 1, 14\}, \{2, 2, 13\}, \{3, 3, 1\}, \{4, 4, 2\}, \{5, 5, 12\}, \\
& \{6, 6, 11\}, \{7, 8, 4\}, \{8, 7, 3\}, \{9, 11, 8\}, \{10, 12, 7\}, \\
& \{11, 9, 6\}, \{12, 10, 5\}, \{13, 13, 10\}, \{14, 14, 9\}],\\
& Z[\{1, 8, 4\}, \{2, 7, 3\}, \{3, 12, 7\}, \{4, 11, 8\}, \{5, 10, 5\}, \\
& \{6, 9, 6\}, \{7, 14, 9\}, \{8, 13, 10\}, \{9, 6, 11\}, \{10, 5, 12\}, \\
& \{11, 2, 13\}, \{12, 1, 14\}, \{13, 4, 2\}, \{14, 3, 1\}]\}\\
$D_{-4}(5_2)$ & \\
\hline
\end{tabular}

\end{center}

\end{document}